\documentclass[10pt,reqno]{amsart}
\usepackage{amsthm,amssymb,amsmath}
\usepackage{graphicx}

\newtheorem{lemma}{Lemma}
\newtheorem{theorem}{Theorem}
\newtheorem{corollary}{Corollary}

\newtheorem{definition}{\rm D\,e\,f\,i\,n\,i\,t\,i\,o\,n}

\theoremstyle{definition}

\newtheorem{example}{Example}

\newcommand{\supp}{\mbox{supp}}
\newcommand{\ordin}{\mbox{\rm{ord}}}

\newcommand{\cl}{\mbox{cl}}

\newcommand{\loc}{\rm loc \it}
\newcommand{\var}{\mbox{var}}

\newcommand{\dyn}{D}
\newcommand{\us}{U}

\begin{document}


\title[Systems with distributions and viability theorem]
{Systems with distributions
and viability theorem}


%
%
%

\author[D.~Kinzebulatov] {D.~Kinzebulatov}

\address{Department of Mathematics and Statistics, University of Calgary \\ 2500 University Drive N.W., Calgary, Alberta, Canada T2N 1N4}

\email{knzbltv@udm.net}


\begin{abstract}
The approach to the consideration of the ordinary differential equations with distributions in the classical space $\mathcal D'$ of distributions with continuous test functions
has certain insufficiencies: the notations are incorrect from the 
point of view of distribution theory, the right-hand side has 
to satisfy the restrictive conditions of equality type. 
In the present paper we consider an initial value problem for the ordinary differential equation with distributions in the space of distributions with dynamic test functions $\mathcal T'$, where the continuous operation of multiplication of distributions by discontinuous functions is defined \cite{DerKin3}, and show that this approach does not have the aforementioned insufficiencies.
We provide the sufficient conditions for viability of solutions of the ordinary differential equations with distributions (a generalization of the Nagumo Theorem), and show that the consideration of the distributional (impulse) controls in the problem of avoidance of encounters with the set (the maximal viability time problem) allows us to provide for the existence of solution, which may not exist for the ordinary controls. 
\end{abstract}

\maketitle



\section{Introduction}
At the end of the nineteenth century K.~Weierstrass pointed out the existence of certain problems in the calculus of variations
where the minimum can not be achieved by 
any smooth function, and the minimizing 
sequence
converges to a discontinuous function (e.g. the Goldschmidt solution of the surface of revolution area minimization problem \cite{Akh}).
As it was further observed, the 
existence of the discontinuous solutions is not a special, but a general situation,
for a large class of problems of optimal control theory (viewed as generalizations of the problems of calculus of variations \cite{Clar}) with one-sided phase constraints on control \cite{Kra, War}. The necessity of extension of 
these problems 
in order to provide the existence of solutions, leads to
ordinary differential equations with distributions or, 
equivalently, 
ordinary differential equations with measures \cite{Art2,Art,Mil,Ses,Silv}. 
\begin{example}[\cite{Dyh}]
\label{ex1} 
The ordinary differential equation with distributions arises in the following optimal control theory problem
\begin{equation}
\label{controlproblem}
\dot{x}=x+v, \quad x(0-)=0, \quad x(1-)=1, \quad v \geq 0,
\end{equation}
\begin{equation}
\label{controlproblem2}
J(v)=\int_{(0,1)} v(t)dt \underset{v}{\to} \min,
\end{equation}
where $x$ is the coordinate of the space shuttle, $v$ is the current charge of fuel,
$J$ determines the total charge of fuel. 
Problem (\ref{controlproblem})(\ref{controlproblem2}) does not have a solution among the pairs $(v,x)$ of locally summable controls $v \in \mathbb L^{\loc}$ 
and locally absolutely-continuous $x \in \mathbb{AC}^{\loc}$, 
and the minimizing sequences 
converge to the elements
\begin{equation}
v^*=\frac{1}{2}\delta_0, 
\quad x^*(t)=\left\{
\begin{array}{ll}
0, &t<0, \\
e^{t-1}, & t>0,
\end{array}
\right.
\end{equation}
where $\delta_0$ is the delta-function concentrated at point $t=0$.
\end{example}

Let us consider the following ordinary differential equation
\begin{equation}
\label{ieq1}
\dot{x}=f(t,x)+g(t,x)v,
\end{equation}
where $v$ is a distribution, which was studied, in particular, in \cite{Art2,Art,Bre,Dyh,Mil,Ses,Silv}, in the classical space $\mathcal D'$ of distributions with continuous test functions \cite{Shi}. 
Since $v \in \mathcal D'$ can be represented as the limit of a sequence of ordinary functions $\{v_k\}_{k=1}^\infty \subset \mathbb L^{\loc}$, a solution of 
(\ref{ieq1}) can be viewed as the limit
\begin{equation}
\label{intro_sol}
x := \lim_{k \to \infty} x_k \text{ in } \mathcal D',
\end{equation}
where $x_k \in \mathbb{AC}^{\loc}$ is the solution of the ordinary differential equation
\begin{equation}
\label{ieq2}
\dot{x}_k=f(t,x_k)+g(t,x_k)v_k(t).
\end{equation}
So any solution of (\ref{ieq1}) is a function of locally bounded variation, in general case discontinuous \cite{Art2,Art,Bre,Dyh,Mil,Ses,Silv}.

However, the consideration of the equation (\ref{ieq1}) in the space $\mathcal D'$ leads to certain problems.

Namely, the notations in (\ref{ieq1}) are incorrect from the point of view of distribution theory, since (\ref{ieq1}) contains the product of a discontinuous function $g\bigl(\cdot,x(\cdot)\bigr)$ and a distribution $v \in \mathcal D'$, which is undefined in $\mathcal D'$. 

Further, let $n$ be the dimension of the system (\ref{ieq1}). If $n \geq 2$, then the solution 
(\ref{intro_sol}) is independent of the choice of $\{v_k\}_{k=1}^\infty$ if and only if the condition of equality type (Frobenius condition) for the right-hand side of (\ref{ieq1}) is satisfied \cite{Mil}. However, this condition can not be satisfied for a large class of problems of optimal control theory, e.g., see \cite{Bre,Mil}.
The approach to consideration of the sequence $\{v_k\}_{k=1}^\infty$ as the part of the system (\ref{ieq1}) is not new (e.g., see \cite{Bre,Mil}),
however,
it was not studied from the point of view of distribution theory, since the choice of a particular sequence $\{v_k\}_{k=1}^\infty$ 
is not motivated by any properties of the space $\mathcal D'$. 

Let us note that since $\mathcal D'$ is isomorphic to the space of Borel measures \cite{Shi} on $\mathbb R$, the problems described above can be reformulated for the ordinary differential equations with measures. 

The insufficiencies described above are due to the absence in $\mathcal D'$ of the continuous operation of multiplication of distributions by discontinuous functions. 
In the present paper we overcome the aforementioned insufficiencies by consideration of the ordinary differential equation with distributions (\ref{ieq1}) in the space $\mathcal T'$ of distributions with \textit{dynamic test functions} \cite{DerKin3}, where the continuous operation of multiplication by discontinuous functions is defined (see the discussion on the problem of multiplication of distributions by discontinuous functions in the classical space $\mathcal D'$ and in other approached, in particular, in the Colombeau generalized functions algebra \cite{Col1,ColM1,Obe} therein).

Also, we study the property of \textit{viability} of solutions of the ordinary differential equations with distributions in $\mathcal T'$. The notion of a viable solution was introduced by J.-P.~Aubin \cite{Aub,Aub2} for the following ordinary differential equation
\begin{equation}
\label{ieqviab}
\dot{x}=f(t,x),
\end{equation}
where $f$ is continuous in $t \in I$ and locally Lipschitz in $x \in \mathbb R^n$.
Let $M \subset \mathbb{R}^n$ be a closed subset, $t_0 \in I$, $\Omega=(t_0,T) \subset I$ is an open interval, in general case unbounded. Following \cite{Aub2}, we give the next definition.

\begin{definition}
\label{viabdef}
A solution of \rm (\ref{ieqviab}) \it such that 
$x(t_0) \in M$ and
\begin{equation}
x(t) \in M 
\end{equation}
for all $t \in \Omega$ is said to be \textit{viable} in $M$ \rm (on $\Omega$). \it The set $M$ is said to have the \textit{property of viability} \rm(\it on $\Omega$\rm) \it for \rm (\ref{ieqviab}), \it if any solution of \rm (\ref{ieqviab}) \it such that $x(t_0) \in M$ is viable in $M$ \rm(\it on $\Omega$\rm). 
\end{definition}

The property of viability is closely related to the problems of existence of the equilibrium points, construction of non-smooth Liapunov functionals, and the problems of optimal control theory \cite{Aub5}, in particular, 
the problem of construction of the admissible control $v_* \in \mathcal V$ for a controlled system
\begin{equation}
\label{intro_contr2}
\dot{x}=f(t,x,v), \quad x(t_0)=x_0, \quad v \in \mathcal V,
\end{equation}
such that for $v=v_*$ the solution of (\ref{intro_contr2}) is viable in $M$ on $\Omega=(t_0,T_*)$, where $T_*=\max_{v \in \mathcal V} T$ (the \textit{problem of avoidance of encounters} with the set $\mathbb{R}^n \setminus M$, or the \textit{maximal viability time problem} \cite{Faz}). 

For ordinary differential equations the sufficient condition for the set $M$ to have the viability property
(the necessary and sufficient condition in the autonomous case) was given in \cite{Nag}. 
\begin{definition}[\cite{Aub2}]
\label{bouligand}
Let $\rho_M:\mathbb R^n \to \mathbb R$ be the distance function, $\rho_M(x) := \inf_{c \in M}\{|x-c|\}$.
The set 
\begin{equation}
K_M(x)=\{y \in {\mathbb R}^n:~ \liminf_{\varepsilon \to 0+} \bigl(\rho_M \bigl(x+\varepsilon y\bigr)/\varepsilon\bigr)=0\}
\end{equation}
is called the contingent cone to the set $M$ at the point $x \in \mathbb{R}^n$.
\end{definition}

\begin{theorem}[Nagumo Theorem]
\label{nagteo2}
Let $M \subset \mathbb{R}^n$ be closed, $M \ne \varnothing$. If 
\begin{equation}
\label{itang}
f(t,x) \in K_M(x)
\end{equation}
for all $t \in \Omega$, $x \in \partial M$, where $\partial M$ is the boundary of $M$, then $M$ has the property of viability for \rm{(}\ref{ieqviab}\rm{)} \it on $\Omega$.
\end{theorem}

If the system (\ref{ieqviab}) is autonomous, then (\ref{itang}) is a necessary and sufficient condition for $M$ to have the property of viability for \rm{(}\ref{ieqviab}\rm{)} (on $\Omega$).

After the ordinary differential equations, the conditions for viability were obtained for the differential inclusions \cite{Aub3,Aub2}, stochastic differential equations \cite{Aub7,Aub4}, differential equations and differential inclusions with aftereffect \cite{Bar,Had2}.
In the present paper we consider the property of viability for the systems with distributions
\begin{equation}
\label{viab_eq5}
\dot{x}=f(t,x)+g(t,x)v, \quad v \in \mathcal T',
\end{equation}
and provide the sufficient condition for a closed set $M$ which is given by the analytic constraints
\begin{equation}
M=\{x \in \mathbb{R}^n: \eta_i(x) \leq 0,~ 1 \leq i \leq m\},
\end{equation}
where $\eta_i:\mathbb R^n \to \mathbb R$ are continuously differentiable functions ($1 \leq i \leq m$), to have the property of viability for the system (\ref{viab_eq5}) (a generalization of Nagumo's Theorem). 
As an application of the results obtained, we provide the sufficient condition for the uniform stability of the equilibrium points of the systems with distributions of the form
\begin{equation*}
\dot{x}=f(x)+g(x)v, \quad v \in \mathcal T'.
\end{equation*}

The consideration of the property of viability for the system (\ref{viab_eq5}) is also motivated by the problem of avoidance of encounters with the set $\mathbb{R}^n \setminus M$. Namely, the problem of avoidance of encounters may have no solution for a system with the ordinary controls
\begin{equation}
\label{viab_eq6}
\dot{x}=f(t,x)+g(t,x)v, \quad v \in \mathbb L,
\end{equation}
which might be viewed as the restriction of the system (\ref{viab_eq5}) to the set of regular controls $v$, and have a solution for the extended system (\ref{viab_eq5}), i.e., when control is allowed to be distributional (impulse). Furthermore, the optimal solution of (\ref{viab_eq5}) has a natural interpretation as the limit of a sequence of solutions of (\ref{viab_eq6}). 
Also, let us mention that the consideration of the property of viability for the system (\ref{viab_eq5}) in the space $\mathcal T'$ allows us to consider the trajectory at the moment of discontinuity (i.e., at the moment of the concentration of the delta-function at $v$), which is important for the definition of the property of viability.

%

\subsection{Notations}

Let $I=(a,b) \subset \mathbb R$ be an open interval, $-\infty \leqslant a<b\leqslant \infty$. We denote by $\mathbb{R}^n$ and $\mathbb{R}^{n \times n}$ the spaces of the vectors and the square matrices of order $n \in \mathbb N$, respectively, with the elements from $\mathbb R$, endowed with the $\max$-norms $|\cdot|$.
Further $(\cdot,\cdot)$ stands for the inner product, $\langle \cdot,\cdot\rangle$ stands for the componentwise product in $\mathbb R^n$. For any $p$, $q \in \mathbb R^n$ we have $|\langle p,q\rangle| \leq |p||q|$.
In what follows, we denote by $\hookrightarrow$ the embedding, i.e., the map preserving the linear and the topological structure.

Let $\hat{\mathbb G}$ be the algebra of bounded functions $g:I \to \mathbb R$ possessing the one-sided limits 
\begin{equation}
g(a+),~g(b-),~g(t+),~g(t-)
\end{equation}
for all $t \in I$.
We identify the elements of $\hat{\mathbb G}$ having the same one-sided limits, and denote the obtained algebra of functions (i.e., the equivalence classes of such functions) by $\mathbb G=\mathbb G(I)$. We define a norm
\begin{equation}
\|g\|_{\mathbb G}:=\sup_{t \in I}\max\{|g(t+)|,|g(t-)|\},
\end{equation}
so $\mathbb G$ is a Banach algebra \cite{Der2}. Following \cite{Dieu}, we call the elements of $\mathbb G$ the \textit{regulated functions}. 
Let 
$T(g):=\{t \in I: \sigma_\tau(g) := g(\tau+)-g(\tau-) \ne 0\}$ be set of points of discontinuity of a regulated function $g \in \mathbb G$.
\begin{lemma}[\cite{Der2}]
Let $g \in \mathbb G$. Then $T(g)$ is at most countable. 
\end{lemma}

Given $g \in \mathbb G$ let us consider a partition $d=\{\tau_i\}_{i=1}^n$, $a<\tau_1<\dots<\tau_n<b$. We define
\begin{equation}
\var_{I}(g):=\sup_{d}\sum_{k=1}^{n-1}|g(\tau_{k+1}-)-g(\tau_k+)|
\end{equation}
(the \textit{total variation}). If $\var_{I}(g)<\infty$, then we call $g$ the \textit{function of bounded variation}. We denote the algebra of functions of bounded variation by $\mathbb{BV}=\mathbb{BV}(I)$. We define
\begin{equation}
\|g\|_{\mathbb{BV}}:=|g(a+)|+\var_{I}(g),
\end{equation}
so $\mathbb{BV}$ is the Banach algebra.

Let $\mathbb C=\mathbb C(I)$ be the subalgebra of the continuous elements of $\mathbb G$. Clearly, $\mathbb C$ consists of the bounded continuous functions on $I$.

Let $\mathbb{CBV}$ be the subalgebra of continuous elements of $\mathbb{BV}$.
If $g \in \mathbb{BV}$, then we denote by $g_c \in \mathbb{CBV}$ the continuous part of $g$ \cite{Der2}.

We denote by $\mathbb L=\mathbb L(I)$ the Banach algebra of functions which are Lebesgue summable on $I$. Let us denote by $\mathbb{AC}=\mathbb{AC}(I)$ the Banach algebra of functions which are absolutely continuous on $I$. 
By $\mathbb{AC}^{\loc}$ and $\mathbb{L}^{\loc}$ we denote the algebras of locally absolutely continuous functions and locally summable functions, respectively \cite{Dun}. We denote by $\mathbb{BV}^{\loc}$ and $\mathbb{CBV}^{\loc}$ the algebras of functions of locally bounded variation.

Let us denote by $\mathbb F=\mathbb F(I)$ the algebra 
of functions $I \to \mathbb R$. We call the elements of $\mathbb F$ the \textit{ordinary functions}.

The definitions given above can be transferred without significant changes to the case of a finite closed interval.

\section{Distributions}

Let us recall the construction of the space of distributions with dynamic test functions \cite{DerKin3}.

\subsection{Dynamic functions} Let $J=\bigl[-\frac{1}{2},\frac{1}{2}\bigr]$.
\begin{definition}[\cite{DerKin3}]
A map
$f:I \to \mathbb F(J)$
is called the \textit{dynamic function}.
\end{definition}
We denote the set of dynamic functions by $d\mathbb F=d\mathbb F(I)$,
and define in $d\mathbb F$ the pointwise operations of addition, multiplication by the element of $\mathbb R$ and multiplication,
so $d \mathbb F$ forms an algebra.

For a given $f \in d\mathbb F$ we call $f(t)(\cdot)$ a \textit{dynamic value} of $f$ at $t$. If $f(t)(\cdot)$ is identically equal to a constant, then we call $f(t)(\cdot)$ an \textit{ordinary value}, and denote it by $f(t)$. 

We define the embedding $\mathbb F \hookrightarrow d\mathbb F$ as follows. For a given $\hat{f} \in \mathbb F$ we associate $f \in d\mathbb F$ such that $f(t)(\cdot) \equiv \hat{f}(t)$ for all $t \in I$. 

Let us denote
\begin{displaymath}
\us(f)=\{t \in I: f(t)(\cdot) \text{ is an ordinary value}\}, \quad \dyn(f)=I \setminus \us(f).
\end{displaymath}

For a given dynamic function $f \in d\mathbb F$ and an ordinary function $g \in \mathbb F$ we define the composition $g \circ f \in d\mathbb F$ by the formula
\begin{equation}
\label{composition}
\bigl(g \circ f\bigr)(t)(\cdot):= g\bigl(f(t)(\cdot)\bigr)
\end{equation}
for all $t \in I$. Accordingly, define an absolute value $|f| \in d\mathbb F$ of $f \in d\mathbb F$ by 
\begin{equation*}
|f|(t)(\cdot):=|f(t)(\cdot)| 
\end{equation*}
for all $t \in I$.
Let us define the support of a dynamic function $f \in d\mathbb F$ by
\begin{equation*}
\supp(f):=\cl\{t \in I: t \in \dyn(f)\text{ or }f(t) \ne \rm 0\} \it \subset I.
\end{equation*}
Let $K \subset I$. We define
\begin{equation*}
\sup_{K}(f) := \sup_{t \in K}\sup_{s \in J}\bigl(f(t)(s)\bigr), \quad \inf_{ K}(f) := \inf_{t \in K}\inf_{s \in J}\bigl(f(t)(s)\bigr).
\end{equation*}

A dynamic function $f \in d\mathbb F$ is said to be bounded on $K \subset I$, if $\sup_{K}|f|<\infty$.
A dynamic function $f \in d\mathbb F$ is said to be non-negative on $K$ {\rm (}denote $f \geq 0$ on $K${\rm )}, if $\inf_K f \geq 0$ (analogously define a non-positive dynamic function).
A constant $c \in \mathbb R$ is called the right-sided limit of $f$ at $\tau \in I$ {\rm (}denote $f(\tau+)=c${\rm )}, if for any $\varepsilon>0$ there exists $\eta>0$ such that 
$\sup_{t \in (\tau,\tau+\eta)}|f-c| \leq \varepsilon$
(analogously define the left-sided limit at $\tau \in I=(a,b)$, and the one-sided limits at points $a, b$).
A dynamic function $f \in d\mathbb F$ is said to be continuous at $\tau \in I$, if $\tau \in \us(f)$ and
$f(\tau)=f(\tau+)=f(\tau-)$
(otherwise we say that $f$ is discontinuous at $\tau \in I$,
denote the set of points of discontinuity of $f$ by $T(f) \subset I$; notice, that $\dyn(f) \subset T(f)$ for any $f \in d\mathbb F$).

Notice that if the dynamic functions in the definitions given above are the ordinary ones, then these definitions coincide with the ordinary definitions. 

For the purpose of construction of the space of distributions, we are mainly interested in the following algebras of dynamic functions.
Suppose that $g \in d\mathbb F$ is such that $g(t)(\cdot) \in \mathbb G(J)$ ($t \in I$), $g$ possesses one-sided limits 
$g(a+)$, $g(b-)$, $g(t+)$, $g(t-)$
for any $t \in I$. The algebra of such dynamic functions (the \textit{regulated dynamic functions}) is denoted by $d\mathbb G$, and endowed with the norm
\begin{equation}
\|g\|_{d\mathbb G} := \sup_{I}|g|,
\end{equation}
so $d\mathbb G$ is the Banach algebra. 
\begin{lemma}[\cite{DerKin3}]
\label{lem1}
Let $g \in d\mathbb G$. Then $T(g)$ is at most countable.
\end{lemma}

Since $\dyn(g) \subset T(g)$, we have that any $g \in d\mathbb G$ possesses the ordinary values everywhere on $I$ except for certain at most countable set.
For a given $g \in d\mathbb G$, we define an ordinary function $\hat{g}$ by
$\hat{g}(t) := g(t)$
($t \in \us(g)$).
We call $\hat{g}$ the \textit{ordinary part} of $g \in d\mathbb G$, and denote 
$\ordin(g) := \hat{g}$.
\begin{lemma}[\cite{DerKin3}] Let $g \in d\mathbb G$. Then
$\hat{g}=\ordin(g)$ is defined everywhere on $I$ for except for at most a countable set, is an element of $\mathbb G$, and
\begin{equation*}
g(t+)=\hat{g}(t+), \quad g(t-)=\hat{g}(t-)
\end{equation*}
for all $t \in I$.
\end{lemma}

Let $g \in d\mathbb G$ be such that
$g(t)(\cdot) \in \mathbb{AC}(J)$, $g(t)(-1/2)=g(t-)$, $g(t)(1/2)=g(t+)$
($t \in I$). We denote the algebra of such dynamic functions by $s\mathbb G$. Then $\mathbb C \hookrightarrow s\mathbb G \hookrightarrow d\mathbb G$. Let us note that for any $g \in s\mathbb G$ we have $\dyn(g)=T(g)$.

Let $g \in s\mathbb G$ be such that $\ordin(g) \in \mathbb{BV}$ and
\begin{equation*}
\sum_{t \in \dyn(g)}\var_{s \in J}\bigl(g(t)(s)\bigr)<\infty.
\end{equation*}
We denote the Banach algebra of such 
dynamic functions (the \textit{dynamic functions of bounded variation}) by $s\mathbb{BV}$, and endow it with the norm
\begin{equation*}
\|g\|_{s\mathbb{BV}}:= |g(a+)|+\|g_c\|_{\mathbb{CBV}}+\sum_{t \in \dyn(g)}\var_{s \in J}\bigl(g(t)(s)\bigr),
\end{equation*}
where $g_c \in \mathbb{CBV}$ is the continuous part of $\ordin(g) \in \mathbb{BV}$.  
We have that $\mathbb{CBV} \hookrightarrow s\mathbb{BV}$.
We define analogously the algebra $s\mathbb{BV}^{\loc}$ of dynamic functions of locally bounded variation.

\begin{example} The Heaviside function $\theta_\tau^\beta \in s\mathbb{BV}$ is defined by 
\begin{equation*}
\theta_\tau^\beta(t)=\left\{
\begin{array}{l}
1, \quad t>\tau, \\
0, \quad t<\tau,
\end{array}
\right. \qquad
\theta_\tau^\beta(\tau)(\cdot)=\beta(\cdot), \\
\end{equation*}
where $\beta \in \mathbb{AC}(J)$ is such that $\beta(-1/2)=0$,
$\beta(1/2)=1$. 
\end{example}

Let $\mathbb G_n$, $\mathbb{AC}_n$, $\mathbb{BV}_n$, $\mathbb L_n$ and $s\mathbb{G}_n$, 
$s\mathbb{BV}_n$ be the spaces of vector-valued functions 
and dynamic functions, respectively, with the operations defined componentwise, and the norms agreeing with the norm in $\mathbb R^n$.

\subsection{Distributions} We denote by $\mathcal D$ the space of the continuous test functions, i.e., the space of the elements of $\mathbb C$ having compact support in $I$, and endowed with the standard locally-convex topology \cite{Shi}.
We denote by $\mathcal T$ the space of the elements $\varphi \in d\mathbb G$ having compact support $\supp(\varphi) \subset I$. The space $\mathcal T$ is called the \textit{space of dynamic test functions} \cite{DerKin3}. 
The sequence $\{\varphi\}_{k=1}^\infty$ is said to be convergent to $\varphi$ in $\mathcal T$, if $\varphi_k \to \varphi$ in $d\mathbb G$ and there exists a compact set $K \subset I$ such that $\supp(\varphi_k) \subset K$ for all $k \in \mathbb N$. Clearly, we have the embedding $\mathcal D \hookrightarrow \mathcal T$. 

We use the notation $\mathcal D'$ for the space of classical distributions, i.e., the space of continuous linear functionals $\mathcal D \to \mathbb R$ \cite{Shi}.
We denote by $\mathcal T$ the space of distributions with dynamic test functions, i.e., the space of continuous linear functionals $\mathcal T \to \mathbb R$ \cite{DerKin3}. The value of a distribution $f \in \mathcal T'$ on a test function $\varphi \in \mathcal T$ is denoted by $(f,\varphi) \in \mathbb R$.
The linear operations in $\mathcal T'$ are defined in a standard way, the space $\mathcal T'$ is endowed with the weak topology, so $f_k \to f$ in $\mathcal T'$ if and only if for any $\varphi \in \mathcal T$ $(f_k,\varphi) \to (f,\varphi)$.

\begin{theorem}[\cite{DerKin3}]
Any distribution in $\mathcal D'$ admits an extension from $\mathcal D$ to $\mathcal T$.
\end{theorem}

\begin{example}
Given $f \in \mathbb L^{\loc}$, let us define a \textit{regular distribution}
by the formula
\begin{equation*}
(f,\varphi) := \int_I f(t)\hat{\varphi}(t)dt,
\end{equation*}
where $\varphi \in \mathcal T$, $\hat{\varphi}=\ordin(\varphi)$. Since $\mathcal D \hookrightarrow \mathcal T$, and the linear manifold of regular distributions in $\mathcal D'$ is isomorphic to $\mathbb L^{\loc}$ \cite{Shi}, we may identify the elements of $\mathbb L^{\loc}$ and the corresponding regular distributions in $\mathcal T'$.
\end{example}

\begin{example}
The \textit{delta-function} $\delta_\tau^\alpha \in \mathcal T'$ is defined by the formula
\begin{equation}
\label{deltafunc}
(\delta_\tau^\alpha,\varphi):=\int_J \alpha(s)\varphi(\tau)(s)ds,
\end{equation}
where $\varphi \in \mathcal T$, $\tau \in I$, and the parameter $\alpha \in \mathbb L(J)$ such that 
\begin{equation}
\label{normcond}
\int_J\alpha(s)ds=1, 
\end{equation}
is called the \textit{shape} of the delta-function. 
Notice, that for any $\varphi \in \mathcal D$ we have
\begin{equation*}
(\delta_\tau^\alpha,\varphi)=\int_J\alpha(s)\varphi(\tau)ds=\varphi(\tau), 
\end{equation*}
so $\delta_\tau^\alpha$ is an extension of the classical delta-function $\delta_\tau \in \mathcal D'$ from $\mathcal D$ to $\mathcal T$.

For a given $\alpha$ satisfying (\ref{normcond}) we define
a sequence $\{\omega_n^\alpha\}_{n=1}^\infty$,
\begin{equation}
\label{deltasequence}
\omega_n^\alpha(t):=\left\{
\begin{array}{ll}
n\alpha(n(t-\tau)), & t \in \bigl(\tau-\frac{1}{2n},\tau+\frac{1}{2n}\bigr), \\
0, &\text{otherwise},
\end{array}
\right. 
\end{equation}
which is called the \textit{delta-sequence having the shape} $\alpha$.
\end{example}

We define the product of $f \in \mathcal T'$ and $g \in d\mathbb G$ by the formula
\begin{equation}
\label{product}
(gf,\varphi) := (f,g\varphi), 
\end{equation}
where $g\varphi \in \mathcal T$. The operation of multiplication defined by (\ref{product}) is continuous, commutative and associative in the sense that $(hg)f=h(gf)$ in $\mathcal T'$ for any $h$, $g \in d\mathbb G$, $f \in \mathcal T'$ \cite{DerKin3}.

\begin{example} 
As follows from (\ref{product}), the product of the Heaviside function $\theta_\tau^\beta \in s\mathbb G$ and the delta-function $\delta_\tau^\alpha \in \mathcal T'$ is defined by the formula
\begin{equation*}
\theta_\tau^\beta \delta_\tau^\alpha=\biggl(\int_J \beta(s)\alpha(s)ds\biggr)\delta_\tau^\gamma,
\end{equation*}
where $\gamma(s)=\alpha(s)\beta(s)/\int_J \beta(s)\alpha(s)ds$ ($s \in J$) satisfies (\ref{normcond}), if $\int_J \beta(s)\alpha(s)ds \ne 0$.
\end{example}

Let $f \in \mathcal T'$. Let us define the support $\supp(f) \subset I$ to be the minimal closed set such that for any $\varphi \in \mathcal T$ with $\supp(\varphi) \cap \supp(f)=\varnothing$ we have $(f,\varphi)=0$.

A distribution $f \in \mathcal T'$ is called non-negative (non-positive) if for any $\varphi \in \mathcal T$ such that $\varphi \geq 0$ we have
$(f,\varphi) \geq 0$ (or $(f,\varphi) \leq 0$, respectively). 

For a given distribution $f \in \mathcal T'$ let us define the value of the integral over $(t_0,t_1)$ ($[t_0,t_1] \subset I$) by the formula
\begin{equation}
\label{integralT}
\int_{(t_0,t_1)} fdt := (f,\varphi_{t_0,t_1}),
\end{equation}
where $\varphi_{t_0,t_1}(t)=1$ if 
$t_0 <t<t_1$, and $\varphi_{t_0,t_1}=0$ otherwise, so $\varphi_{t_0,t_1} \in \mathcal T$.
The operator of integration in $\mathcal T'$ which is given by (\ref{integralT}) is linear and continuous \cite{DerKin3}. 

Let
$g \in s\mathbb{BV}^{\loc}$, we define the derivative $\dot{g} \in \mathcal T'$ by the formula
\begin{equation}
\label{deriv}
(\dot{g},\varphi):= \int_I \hat{\varphi}(t)dg_c(t)+\sum_{\tau \in T(g)}\int_J\varphi(\tau)(s)(g(\tau)(s))^{\cdot}_sds.
\end{equation}
where $\varphi \in \mathcal T$, $g_c \in \mathbb{CBV}^{\loc}$, and the set of points of discontinuity $T(g)$ is at most countable by Lemma \ref{lem1}.
If for any $\tau \in T(g)$ $\sigma_\tau(g)=g(\tau+)-g(\tau-) \ne 0$, then
\begin{equation}
\label{deriv2}
\dot{g}=\dot{g}_c+\sum_{\tau \in T(g)}\sigma_{\tau}(g)\delta_{\tau}^{\alpha_\tau}
\end{equation}
in $\mathcal T'$, where the shape of a delta-function
$\alpha_\tau(s)=(g(\tau)(s))^{\cdot}_s/\sigma_{\tau}(g)$ 
($s \in J$).

\begin{example}
The derivative of the Heaviside function $\theta_\tau^\beta$ is the delta-function $\delta_\tau^\alpha$,
\begin{equation*}
\dot{\theta}_\tau^\beta=\delta_\tau^\alpha,
\end{equation*}
where $\alpha=\dot{\beta}$.
\end{example}

In what follows the notations $\mathcal D_n'$ and $\mathcal T_n'$ stand for the spaces of the vector-valued distributions with the convergence and the operations defined componentwise.

\section{Systems with distributions \\ [3mm]} 
Let $D \subset I \times \mathbb R^n$ be open.
Consider in $\mathcal T_n'=\mathcal T_n'(I)$ an initial value problem
\begin{equation}
\label{eq1}
\dot{x}=f(t,x)+g(t,x)v, \quad x(t_0-)=x_0,
\end{equation}
where $(t_0,x_0) \in D$, 

1) The function $f:D \to \mathbb{R}^n$ is Caratheodory in $(t,x) \in D$ and Lipschitz in $x$ with the constant $K_f>0$ in $D$.

2) The function $g:D \to \mathbb{R}^{n \times n}$ belongs to $s\mathbb G_n$ in $t$ (in particular, $g$ is continuous in $t$) and Lipschitz in $x$ with the constant $K_g>0$ in $D$.

3) The distribution $v \in \mathcal T_n'$ is defined by
\begin{equation*}
v=\dot{u}, \\ [2mm]
\end{equation*}
where $u \in s\mathbb{BV}_n^{\loc}$, we assume that for any $\tau \in T(u)$ $\sigma_\tau(u) \ne 0$, so according to the definition of the derivative in $\mathcal T_n'$
\begin{equation}
\label{vrepr}
v=\dot{u}_c+\sum\limits_{\tau \in T(u)}\bigl\langle 
\sigma_{\tau}(u),\delta_{\tau}^{\alpha_\tau}\bigr\rangle,
\end{equation}
where $\delta_{\tau}^{\alpha_\tau} \in \mathcal T_n'$ is the vector-valued delta-function, 
$(\delta_{\tau}^{\alpha_\tau},\varphi):= \bigl((\delta_{\tau}^{\alpha_\tau^i},\varphi) \bigr)_{i=1}^n$, 
where $\varphi \in \mathcal T$, $\alpha_\tau=(\alpha_\tau^i)_{i=1}^n \in \mathbb L_n(J)$, $\langle,\rangle$ is the componentwise product in $\mathbb{R}^n$.

Let us note that $g$ may be discontinuous in $t$. This is important, since in the problems of optimal control $g$ may also depend on the (discontinuous) ordinary control (in $t$), which is not considered here.

A \textit{solution} of the initial value problem (\ref{eq1}) on $\Omega \subset I$ is the dynamic function $x \in s\mathbb{BV}_n^{\loc}(\Omega)$ such that $x(t)(s) \in D$ ($t \in \Omega$, $s \in J$) and (\ref{eq1}) is satisfied in $\mathcal T_n'(\Omega)$. 
An ordinary part $\hat{x} \in \mathbb{BV}_n^{\loc}(\Omega)$ of the solution $x$ is called the \textit{ordinary solution} of (\ref{eq1}).
Let us note that in contrast to the approach based on the space of distributions $\mathcal D'$ the operations of multiplication, differentiation and composition arising in (\ref{eq1}) are correctly defined in the sense of the distribution theory in $\mathcal T_n'(\Omega)$ (see the Introduction).

In the formulation of the next theorem we put $\Omega=I$.

\begin{theorem}
\label{teo1}
Let $x \in s\mathbb{BV}_n^{\loc}$ be the solution of \rm{(}\ref{eq1}\rm{)}\it, $\hat{x} \in \mathbb{BV}_n^{\loc}$ be the ordinary solution of \rm{(}\ref{eq1}\rm{)}\it. Then
\begin{multline}
\label{inteq}
\hat{x}(t)=x_0+\int_{t_0}^t f\bigl(r,\hat{x}(r)\bigr)dr+\int_{t_0}^t \hat{g}\bigl(r,\hat{x}(r)\bigr)du_c(r)+\\
\sum\limits_{\tau<t} \bigl(\gamma_\tau(1/2)-\hat{x}(\tau-)\bigr)-\sum_{\tau<t_0}(\gamma_\tau(1/2)-\hat{x}(\tau-)),
\end{multline}
and the dynamic value $\gamma_\tau(\cdot) := x(\tau)(\cdot)$ satisfies
\begin{equation}
\label{limitsystem}
\dot{\gamma}_\tau(s)=g\bigl(\tau,\gamma_\tau(s)\bigr)(s)\bigl\langle \sigma_{\tau}(u),\alpha_\tau(s)\bigr\rangle, \quad \gamma_\tau(-1/2)=x(\tau-), 
\end{equation}
where $T(u)=\{\tau\}$ is at most countable by Lemma \ref{lem1}, $g\bigl(t,x\bigr)(s)$ is a dynamic value of $g$ at $(t,x)$, $x(\tau+)=\hat{x}(\tau+)=\gamma_\tau(1/2)$.

Conversely, any $x \in s\mathbb{BV}_n^{\loc}$ satisfying \rm{(}\ref{inteq}\rm{)},\rm{(}\ref{limitsystem}\rm{)} \it is a solution of \rm{(}\ref{eq1}\rm{)}.
\end{theorem}
\begin{proof}
Observe that by our assumptions a solution of the problem (\ref{limitsystem}) exists and is unique. Let us show that the second statement holds.                  
Suppose that there exists a solution $\hat{x} \in \mathbb{BV}_n^{\loc}$ of the integral equation (\ref{inteq}). Then we may define $x \in s\mathbb{BV}^{\loc}$ as in the formulation of the theorem.
First, show that $x \in s\mathbb{BV}_n^{\loc}$ satisfies the initial condition in (\ref{eq1}). We have that
$x(t_0-)=\hat{x}(t_0-)$, the limit $\hat{x}(t_0-)$ exists since $\hat{x} \in \mathbb{BV}_n^{\loc}$. Since all the integrals in (\ref{inteq}) are continuous in $t$, we have
\begin{equation*}
\lim_{t \to t_0-}\int_{t_0}^t f(r,\hat{x}(r))dr=0, \quad \lim_{t \to t_0-}\int_{t_0}^t \hat{g}(r,\hat{x}(r))du_c(r)=0,
\end{equation*}
so the limit of the step part in (\ref{inteq}) is
\begin{equation*}
\lim_{t \to t_0-}\sum_{\tau<t} \bigl(\gamma_\tau(1/2)-\hat{x}(\tau-)\bigr)=\sum_{\tau<t_0}(\gamma_\tau(1/2)-\hat{x}(\tau-)),
\end{equation*}
i.e., $x(t_0-)=x_0$. Second, we show that $x \in s\mathbb{BV}_n^{\loc}$ satisfies the differential equation in (\ref{eq1}) in $\mathcal T'_n$. According to the definition of the derivative, we have
\begin{equation}
(\dot{x},\varphi)=\int_I \hat{\varphi}(t)dx_c(t)+\sum_{\tau \in T(u)}\int_J\varphi(\tau)(s)(x(\tau)(s))^{\cdot}_sds,
\end{equation}
where $x(\tau)(s)=\gamma_\tau(s)$, so differential equation (\ref{eq1}) is equivalent to
\begin{multline}
\notag
\int_I\hat{\varphi}(t)dx_c(t)+\sum_{\tau \in T(u)}\int_J\varphi(\tau)(s)(x(\tau)(s))^{\cdot}_sds=\\=\int_If(t,\hat{x}(t))\hat{\varphi}(t)dt+\int_I\hat{g}(t,\hat{x}(t))du_c(t)+
\sum_{\tau \in T(u)}\int_J\varphi(\tau)(s)g\bigl(\tau,x(\tau)(s)\bigr)(s)\langle \sigma_{\tau}(u),\alpha_\tau(s)\rangle ds.
\end{multline}
Let us show that the equality
\begin{equation}
\label{equality1}
\int_J\varphi(\tau)(s)(x(\tau)(s))^{\cdot}_sds=\int_J\varphi(\tau)(s)g\bigl(\tau,x(\tau)(s)\bigr)(s)\langle \sigma_{\tau}(u),\alpha_\tau(s)\rangle ds
\end{equation}
holds for any $\varphi \in \mathcal T$, $\tau \in T(u)$, and
\begin{equation}
\label{equality2}
\dot{x}_c=f(t,x)+g(t,x)\dot{u}_c.
\end{equation}
Indeed, the equality (\ref{equality1}) follows from (\ref{limitsystem}), where
$x(\tau)(s)=\gamma_\tau(s)$ ($\tau \in T(u)$).
Further, by the definition of the regular distribution in $\mathcal T_n'$, the equality (\ref{equality2}) is equivalent to the equality
\begin{equation*}
\int_I \hat{\varphi}(t)dx_c(t)=\int_I \hat{\varphi}(t)f(t,\hat{x}(t))dt+\int_I \hat{\varphi}(t)\hat{g}(t,\hat{x}(t))du_c(t)
\end{equation*}
for any $\varphi \in \mathcal T_n$. 
Then the last equality is equivalent to
\begin{equation}
\label{lasteq2}
x_c(t)=x_c(t_0)+\int_{t_0}^t f(r,\hat{x}(r))dr+\int_{t_0}^t \hat{g}(r,\hat{x}(r))du_c(r)
\end{equation}
for all $t \in I$. As follows from the construction of $x \in s\mathbb{BV}^{\loc}$ in the statement of the theorem, the value of the jump of the solution $x(\tau+)-x(\tau-)$ can be found from the dynamic value $x(\tau)(\cdot)=\gamma_\tau(\cdot)$. Consequently, the equality of the continuous parts
(\ref{lasteq2}) implies (\ref{inteq}). Since the equalities obtained are equivalent to the equalities (\ref{equality1}),(\ref{equality2}) in
$\mathcal T_n'$, the dynamic function $x \in s\mathbb{BV}_n^{\loc}$ is a solution of the initial value problem (\ref{eq1}). 

We now show that the first statement holds. Let $x \in s\mathbb{BV}_n^{\loc}$ be the solution of the problem (\ref{eq1}). Then (\ref{equality1}) and (\ref{equality2}) are true. Further, the equality (\ref{equality2}) is equivalent to
\begin{equation*}
\int_I \hat{\varphi}(t)dx_c(t)=\int_I \hat{\varphi}(t)f(t,\hat{x}(t))dt+\int_I \hat{\varphi}(t)\hat{g}(t,\hat{x}(t))du_c(t)
\end{equation*}
for all $\varphi \in \mathcal T_n$. 
According to DuBois-Reymond Lemma \cite{Cour}, we have
\begin{equation}
\label{lasteq}
x_c(t)=x_c(t_0)+\int_{t_0}^t f(r,\hat{x}(r))dr+\int_{t_0}^t \hat{g}(r,\hat{x}(r))du_c(r)
\end{equation}
for all $t \in I$. 
Analogously, as follows from DuBois-Reymond Lemma, the equality (\ref{equality1}) implies that (\ref{limitsystem}) holds for any $\tau \in T(u)$. Then (\ref{inteq}) is true, as follows from the definition of the space $s\mathbb{BV}_n^{\loc}$ and the initial condition $x(t_0-)=x_0$.
\end{proof}

\begin{theorem}
\label{existsteo2}
There is a constant $h>0$ such that there exists a solution $x \in s\mathbb{BV}_n(t_0-h,t_0+h)$ of the problem \rm{(}\ref{eq1}\rm{)} \it, which is unique in the sense that it coincides with any other solution of \rm{(}\ref{eq1}\rm{)} \it on the common interval of definition in $(t_0-h,t_0+h)$, and depends continuously on 
$u|_{(t_0-h,t_0+h)} \in s\mathbb{BV}_n(t_0-h,t_0+h)$.
\end{theorem}
\begin{proof}
Let us denote $I_h:=(t_0-h,t_0+h)$, $I_h^+:=(t_0,t_0+h)$, $I_h^-:=(t_0-h,t_0) \subset I=(a,b)$. 
Let us consider in $\mathbb{BV}_n(I_h^+)$ the following integral equation
\begin{equation}
\label{p1}
\hat{x}(t)=x_0+\int_{t_0}^tf\bigl(r,\hat{x}(r)\bigr)dr+\int_{t_0}^t \hat{g}\bigl(r,\hat{x}(r)\bigr)du_c(r)+\sum_{t_0<\tau<t}\bigl(\gamma_\tau(1/2)-\hat{x}(\tau-)\bigr),
\end{equation}
\begin{equation}
\label{p2}
\dot{\gamma}_\tau(s)=g\bigl(\tau,\gamma_\tau(s)\bigr)(s)\bigl\langle\sigma_{\tau}(u),\alpha_\tau(s) \bigr\rangle, \quad \gamma_\tau(-1/2)=\hat{x}(\tau-), \\ [2mm]
\end{equation}
where $\tau \in T(u) \subset I_h^+$, the dynamic function $u$ is assumed to be restricted to $I_h^+$. Clearly, if the solution $\hat{x} \in \mathbb{BV}_n(I_h^+)$ of (\ref{p1})(\ref{p2}) exists, then $\hat{x}(t_0+)=x_0$. Let us show the existence of solution of (\ref{p1})(\ref{p2}) in $\mathbb{BV}_n(I_h^+)$ for some $h>0$. Let $N>0$. We define
\begin{equation}
M_g^h=\max\{|g(t,x)(s)|: s \in J, ~t \in [t_0,t_0+h], ~|x-x_0| \leq N\} \geq 0,
\end{equation}
\begin{equation}
Q^+_h=\{\hat{x} \in \mathbb{BV}_n(I_h^+):\|\hat{x}-x_0\|_{\mathbb{BV}_n(I_h^+)}\leq N\}, 
\end{equation}
so that $[t_0,t_0+h] \times \{x \in \mathbb R^n:|x-x_0| \leq N\} \subset D$. Obviously, $Q^+_h$ is a complete metric subspace. We define a map $P$ on $Q^+_h$ by the formula
\begin{equation*}
P(\hat{x})(t)=x_0+\int_{t_0}^tf\bigl(r,\hat{x}(r)\bigr)dr+\int_{t_0}^t \hat{g}\bigl(r,\hat{x}(r)\bigr)du_c(r)+\sum_{t_0<\tau<t}\bigl(\gamma_\tau(1/2)-\hat{x}(\tau-)\bigr),
\end{equation*}
where $t \in I_h^+$, $\hat{x} \in Q^+_h$. Then
\begin{equation*}
\|P(\hat{x})-x_0\|_{\mathbb{BV}_n(I_h^+)} \leq \int_{I_h^+}\bigl|f\bigl(r,\hat{x}(r)\bigr)\bigr|dr+M_g^h\bigl(\var_{I_h^+}(u)\bigr)^2+M_g^h \var_{I_h^+}(u). 
\end{equation*}
Then there exists $h>0$ such that $\|P(\hat{x})-x_0\|_{\mathbb{BV}_n(I_h^+)} \leqslant N$, i.e., $P:Q^+_h \to Q^+_h$.
Let us show that there exists $h>0$ sufficiently small such that there is $\lambda>0$ such that
\begin{equation}
\label{proof_ineq}
\|P(\hat{x})-P(\hat{y})\|_{\mathbb{BV}_n(I_h^+)} \leq 
\lambda \|\hat{x}-\hat{y}\|_{\mathbb{BV}_n(I_h^+)}
\end{equation}
for all $\hat{x}$, $\hat{y} \in Q^+_h$. Indeed, to have the inequality (\ref{proof_ineq}) satisfied it suffices to put 
\begin{equation}
\lambda=K_fh+K_g\var_{I_h^+}(u)h+K_g\var_{I_h^+}(u)\var_{t_0}^t(u) \geq 0, 
\end{equation}
so we can find $h>0$ sufficiently small such that $\lambda<1$. Then according to the Fixed Point Theorem \cite{Kan} the mapping $P$ has the only fixed point in $Q^+_h$, which is the only solution of the integral equation (\ref{p1})(\ref{p2}). Since the solution $\gamma_\tau$ of the problem (\ref{limitsystem}) exists on $J$ for any $\tau \in T(u)$, and $T(u)=T(\hat{x})=\{t \in I_h^+:\hat{x}(t+) \ne \hat{x}(t-)\}$, we may define $x \in s\mathbb{BV}_n(I_h)$ by $x(t)=\hat{x}(t)$ if $t \not\in T(u)$, and $x(t)(s)=\gamma_t(s)$ ($s \in J$) if $t \in T(u)$. According to Theorem \ref{teo1}, the dynamic function $x \in s\mathbb{BV}_n(I_h^+)$ is the solution of the initial value problem (\ref{eq1}) on $I_h^+$.

Further, by the change of independent variable in (\ref{p1}) from $t$ to $-t$ we can show that there exists $h>0$ such that $[t_0-h,t_0+h] \times \{x \in \mathbb R^n:|x-x_0| \leq N\} \subset D$ and the following integral equation
\begin{multline}
\label{p4}
\hat{x}(t)=x_0+\int_{t_0}^tf\bigl(r,\hat{x}(r)\bigr)dr+\int_{t_0}^t \hat{g}\bigl(r,\hat{x}(r)\bigr)du_c(r)+\\
\sum_{t_0-h<\tau<t}\bigl(\gamma_\tau(1/2)-\hat{x}(\tau-)\bigr)-\sum_{\tau<t_0}\bigl(\gamma_\tau(1/2)-x(\tau-)\bigr),
\end{multline}
endowed with condition (\ref{p2}), where $\tau \in T(u) \subset I_h$, and the dynamic function $u$ is assumed to be restricted to $I_h$, has a unique solution in 
\begin{equation*}
Q_h=\{\hat{x} \in \mathbb{BV}_n(I_h):\|\hat{x}-x_0\|_{\mathbb{BV}_n(I_h)}\leq N\}. 
\end{equation*}

Now let us show that the mapping $P=P_u$ depends continuously on $u \in s\mathbb{BV}_n(I_h^+)$. Suppose that we have $\{u_k\}_{k=1}^\infty \subset s\mathbb{BV}_n(I_h^+)$, 
$u \in s\mathbb{BV}_n(I_h^+)$ and 
\begin{equation*}
u_k \to u
\end{equation*}
in $s\mathbb{BV}_n(I_h^+)$. Let us show that for any $\hat{x} \in Q^+_h$ we have $P_{u_k}(\hat{x}) \to P_u(\hat{x})$ in $Q^+_h$.
Indeed, we have
\begin{multline}
\notag 
\|P_u(\hat{x})-P_{u_k}(\hat{x})\|_{\mathbb{BV}(I_h^+)} 
\leq M_g^h \var_{I_h^+}(u-u_k)
+\\M_g^h\var_{I_h^+}(u-u_k)\var_{I_h^+}(u)
+M_g^h\var_{I_h^+}(u_k)\var_{I_h^+}(u-u_k),
\end{multline}
where $\{\var_{I_h^+}(u_k)\}_{k=1}^\infty$ is bounded, so
\begin{equation*}
\|P_u(\hat{x})-P_{u_k}(\hat{x})\|_{\mathbb{BV}_n(I_h^+)} \to 0 \text{ as }
\|u-u_k\|_{s\mathbb{BV}_n(I_h^+)} \to 0. 
\end{equation*}
As is shown in \cite{Kan}, the inequality above implies the convergence $\hat{x}_k \to \hat{x}$ in $Q^+_h$, where $\hat{x}_k \in Q^+_h$ is a fixed point of the mapping
$P_{u_k}$ ($k \in \mathbb N$). Consequently, the solution $\hat{x} \in \mathbb{BV}_n(I_h^+)$ of the integral equation (\ref{inteq}) depends continuously on $u \in s\mathbb{BV}_n(I_h^+)$. 

Let us show that the solution $x \in s\mathbb{BV}_n(I_h^+)$ of the problem (\ref{eq1})
depends continuously on $u \in s\mathbb{BV}_n(I_h^+)$.
Let us denote the dynamic values $\beta_\tau(\cdot)=u(\tau)(\cdot)$, 
$\beta^k_\tau(\cdot)=u_k(\tau)(\cdot) \in \mathbb{AC}(J)$ ($\tau \in T(u)$), and denote the solution of the system (\ref{limitsystem}) corresponding to $u_k$ by $\gamma_\tau^k \in \mathbb{AC}(J)$.
Employing the known estimations we have
\begin{multline}
\notag
\var_J(\gamma_\tau-\gamma_\tau^k) \leq
M_g^h\biggl(|\sigma_{\tau}(u)-\sigma_{\tau}(u_k)|\var_J(\beta_\tau)+
|\sigma_{\tau}(u_k)|\var_J(\beta_\tau-\beta_\tau^k)\biggr)+\\
K_g|\sigma_{\tau}(u_k)|e^{n\sqrt{n}K_gh}M_g^h\biggl(\var_{I_h^+}(u)|\sigma_{\tau} (u)-\sigma_{\tau}(u_k)|+
|\sigma_{\tau}(u_k)|\var_{I_h^+}(u_k)\var_J(\beta_\tau-\beta_\tau^k)\biggr),
\end{multline}
where the sequences $\{\var_{I_h^+}(u_k)\}$, $\{|\sigma_\tau(u_k)|\}$ are bounded.
Consequently, 
\begin{equation*}
\|x_k-x\|_{s\mathbb{BV}_n(I_h^+)} \to 0 \text{ if } \|u_k \to 
u\|_{s\mathbb{BV}_n(I_h^+)} \to 0. \\ [2mm]
\end{equation*}
The case of the left-sided neighborhood $I_h^-$ is treated similarly, so the solution $x \in s\mathbb{BV}_n(I_h)$ of the initial value problem (\ref{eq1}) depends continuously on $u \in s\mathbb{BV}_n(I_h)$.
\end{proof}

The existence of a non-continuable in $D$ solution of the initial value problem (\ref{eq1}) can be shown in a standard way. \\


Let $D=I \times N$, $N \subset \mathbb R^n$ be open.
An initial value problem
\begin{equation}
\label{eq2}
\dot{x}=f(t,x)+g(t,x)v, \quad x(t_0-)=x_0
\end{equation}
where $v=\dot{u} \in \mathcal D_n'$, $u \in \mathbb{BV}_n^{\loc}$, was considered, in particular, in \cite{Art2,Art,Dyh,Mil,Ses,Silv}.
The solution $x \in \mathbb{BV}_n^{\loc}$ of the initial value problem (\ref{eq2}) is defined by 
\begin{equation*}
x:= \lim_{k \to \infty} x_k \text{  in } \mathcal D_n',
\end{equation*}
where $x_k \in \mathbb{AC}^{\loc}$ is the solution of the initial value problem for the differential equation with the ordinary right-hand side,
\begin{equation*}
\dot{x}_k=f(t,x_k)+g(t,x_k)v_k(t), \quad x_k(t_0)=x_0,
\end{equation*}
where $v_k=\dot{u}_k \in \mathbb L_n^{\loc}$, $u_k \in \mathbb{AC}_n^{\loc}$, and $v_k \to v$ in $\mathcal D_n'$ \cite{Mil,Ses}.
If the function $g$ is continuously differentiable in $(t,x) \in D$, then the necessary and sufficient condition for the uniqueness of the solution of (\ref{eq1}), i.e., its independence on the choice of the approximating sequence $\{v_k\}_{k=1}^\infty$, is the Frobenius condition \\ [-2mm]
\begin{equation}
\label{frob}
\bigl[g^m,g^l\bigr]_x(t,x)=0, \quad (t,x) \in D,~1 \leq l,m \leq n, \\ [1mm]
\end{equation}
where $g^m$ is an $m$-th column of $g$, and $[\cdot,\cdot]_x$ is the Lie bracket in the variable $x$, e.g., see \cite{Mil,Ses}.
If the condition (\ref{frob}) is satisfied, then the solution $x \in \mathbb{BV}_n^{\loc}$ of the initial value problem (\ref{eq2}) satisfies
\begin{multline}
\label{eqses}
x(t)=x_0+\int_{t_0}^t f\bigl(r,x(r)\bigr)dr+\int_{t_0}^t g\bigl(r,x(r)\bigr)du_c(r)+\\
\sum\limits_{\tau<t} \bigl(\gamma_\tau(1/2)-x(\tau-)\bigr)-\sum_{\tau<t_0}\bigl(\gamma_\tau(1/2)-x(\tau-)\bigr),
\end{multline}
\begin{equation}
\label{limitsystem2}
\dot{\gamma}_\tau(s)=g\bigl(\tau,\gamma_\tau(s)\bigr)\sigma_{\tau}(u), \quad \gamma_\tau(-1/2)=x(\tau-), \\ [2mm]
\end{equation}
where $T(u)=\{\tau\}$, $x(\tau+)=\gamma(1/2)$ \cite{Ses}.

As is well known, the condition (\ref{frob}) is a necessary and sufficient condition for the invariance of the value $\gamma(1/2) \in \mathbb{R}^n$ with respect to the choice of the function $\alpha \in \mathbb L_n(J)$ satisfying (\ref{normcond}), where $\gamma$ is a solution of the problem
\begin{equation}
\label{limsyscomplete}
\dot{\gamma}(s)=g\bigl(\tau,\gamma(s)\bigr)\alpha(s), \quad \gamma(-1/2)=\gamma_0,
\end{equation}
where $\tau \in I$, $\gamma_0 \in N$ \cite{Gaishun}. 
Consequently, (\ref{frob}) is equivalent to the condition of the independence of the ordinary solution of problem (\ref{eq1}) on the choice of the shapes of delta-functions $\alpha_\tau$ ($\tau \in T(u)$) (see the Introduction). 

Let us note that if the condition (\ref{frob}) is satisfied, then (\ref{limitsystem2}) coincides with (\ref{limitsystem}) for $\alpha_k \equiv (1,\dots,1)^{\top}$.

\begin{example} Let $I=(-1,1)$, $D=I \times \mathbb R$.
Let us consider in $\mathcal T'$ the initial value problem,
\begin{equation}
\label{eq3}
\dot{x}=\delta_0^\alpha x, \quad x(-1/2-)=x_0,
\end{equation}
where $x_0 \in \mathbb R$. According to Theorem \ref{teo1} the solution $x \in s\mathbb{BV}$ of the problem (\ref{eq3}) is given by
\begin{equation*}
x(t)=\left\{
\begin{array}{ll}
x_0, & t<0, \\
x_0e, & t>0,
\end{array}
\right. \quad
x(0)(s)=x_0e^{\int_{-1/2}^s\alpha(\eta)d\eta} \quad (s \in J), \\ [2mm]
\end{equation*}
where $J=[-1/2,1/2]$.

Note that the same value of the jump of the solution at $t=0$ can be obtained if the delta-function $\delta_0^\alpha$ in (\ref{eq3}) is replaced by the terms of the delta-sequence (\ref{deltasequence}) having the shape $\alpha$. Also note that the value of the jump of the solution at $t=0$ is independent of the choice of the shape $\alpha$ since for $n=1$ the Frobenius condition (\ref{frob}) is always satisfied.
\end{example}

\section{Viability Theorem and stability analysis}


Let $\Omega=(t_0,T) \subset I$, let $M \subset \mathbb R^n$ be a closed subset, $D=I \times N$, where $N$ is an open subset,
$M \subset N$.
The following extends Definition \ref{viabdef}.

\begin{definition}
\label{viabdef2}
A solution of the system \rm (\ref{eq1}) \it
such that $x(t_0-) \in M$ and
\begin{equation*}
x(t)(s) \in M
\end{equation*}
for all $t \in \Omega \cup \{t_0\}$, $s \in J$, where $t_0 \in I$, is said to be viable in $M$ on $\Omega$. The set $M$ is said to have the property of viability for \rm (\ref{eq1}) on $\Omega$, if any solution of \rm (\ref{eq1}) \it such that $x(t_0-) \in M$ is viable in $M$ \rm(\it on $\Omega$\rm).
\end{definition}

Suppose that the set $M \subset \mathbb{R}^n$ is given by
\begin{equation}
\label{M}
M=\{x \in \mathbb{R}^n:\eta_i(x) \leq 0,~ 1 \leq i \leq m\},
\end{equation}
where $\eta_i:\mathbb{R}^n \to \mathbb{R}$ are continuously differentiable on $\mathbb R^n$ ($1 \leq i \leq m$), so $M$ is closed in $\mathbb{R}^n$. 
Clearly, for each $1 \leq i \leq m$ if $x \in \mathbb R^n$ is such that $\eta_i(x)=0$, $\dot{\eta}_i(x) \ne 0$, then 
\begin{equation*}
K_{\{p~:~\eta_i(p)\leq 0\}}(x)=\{y:(\dot{\eta}_i(x),y)\leq 0\}. 
\end{equation*}
Then according to \cite[p.224]{Aub2}
for any $x \in \partial M$ such that $\dot{\eta}_i(x) \in \mathbb{R}^n$ ($i \in L_x=\{i:\eta_i(x)=0\}$) are linearly independent, the contingent cone $K_M(x)$ is given by
\begin{equation}
K_M(x)=\{y \in \mathbb{R}^n:\bigl(\dot{\eta}_i(x),y\bigr)\leq 0,~i \in L_x\}.
\end{equation}
The following theorem follows immediately from the Nagumo Theorem.

\begin{theorem}
\label{nagteo}
Let $\dot{\eta}_i(x) \in \mathbb{R}^n$ \rm(\it $i \in L_x$\rm) \it be linearly independent for any $x \in \partial M$.
If
\begin{equation*}
\bigl(\dot{\eta}_i(x),f(t,x)\bigr) \leq 0 \quad (i \in L_x)
\end{equation*}
for all $t \in \Omega$, $x \in \partial M$, then
$M$ has the the property of viability for \rm(\ref{ieqviab}) \rm(\it on $\Omega$\rm).
\end{theorem}

Let us consider in $\mathcal T_n'$ the following differential equation with distributions
of the form (\ref{eq1}),
\begin{equation}
\label{eq11}
\dot{x}=f(t,x)+g(t,x)v, \quad v=w+\sum\limits_{k=1}^\infty \langle c_k,\delta_{\tau_k}^{\alpha_k} \rangle \in \mathcal T_n',
\end{equation}
where the function $w:I \to \mathbb{R}^n$ is continuous, $\{\tau_k\}_{k=1}^\infty \subset I$, 
$\delta_{\tau_k}^{\alpha_k} \in \mathcal T_n'$, the coefficients $c_k \in \mathbb R^n$ are such that the primitive of $v$ in $\mathcal T_n'$ 
is in $s\mathbb{BV}_n^{\loc}$,
$\alpha_k \in \mathbb C_n(J)$, the functions $f$, $g$ satisfy 1), 2) and $f$ is also continuous in $t \in \Omega$. We suppose that for any $x_0 \in D$ there exists a non-continuable solution of the initial value problem for (\ref{eq11}) in $D$ with the initial value $x(t_0-)=x_0$.

The following statement generalizes Theorem \ref{nagteo}.

\begin{theorem}
\label{viabteo}
Let $\dot{\eta}_i(x) \in \mathbb{R}^n$ \rm(\it $i \in L_x$\rm) \it be linearly independent for any $x \in \partial M$.
If
\begin{equation}
\label{vcond}
\bigl(\dot{\eta}_i(x),f(t,x)+g(t,x)v\bigr)\leq 0 \quad (i \in L_x)
\end{equation}
in $\mathcal T'(\Omega)$ for all $x \in \partial M$, then
$M$ has the property of viability for \rm (\ref{eq11}) \rm(\it on $\Omega$\rm).
\end{theorem}

The definition of a non-positive distribution in $\mathcal T'$ was given in Section 2.
As follows from (\ref{vcond}) and the examples below, the property of viability depends on the choice of the shapes of the delta-function in $v$, including the case where the Frobenius condition (\ref{frob}) is satisfied.

\begin{proof}
1) Let $x \in \partial M$, $i \in L_x$. Let us show that if the inequality (\ref{vcond}) holds, then
\begin{equation}
\label{v1i}
\bigl(\dot{\eta}_i(x),f(t,x)+\hat{g}(t,x)w(t)\bigr)\leq 0
\end{equation}
for all $t \in \Omega$, and
\begin{equation}
\label{v2i}
\bigl(\dot{\eta}_i(x),g(\tau_k,x)(s)\langle c_k, \alpha_k(s)\rangle\bigr)\leq 0
\end{equation}
for all $s \in J$, $k \in \mathbb N$. According to the definition of a non-positive distribution in $\mathcal T'$, the inequality (\ref{vcond}) implies that
\begin{equation*}
\int_I(\dot{\eta}_i(x),f(t,x)+\hat{g}(t,x)w(t))\hat{\varphi}(t)dt+\sum_{k=1}^\infty (\dot{\eta}_i,g(\tau_k,x)(s)\langle c_k,(\delta_{\tau_k}^{\alpha_k},\varphi)\rangle) \leq 0
\end{equation*}
for any $\varphi \in \mathcal T$, $\varphi \geq 0$, where $\hat{\varphi}=\ordin(\varphi)$, i.e., by the definition of the delta-function
\begin{equation}
\label{vineq}
\int_I(\dot{\eta}_i(x),f(t,x)+\hat{g}(t,x)w(t))\hat{\varphi}(t)dt+\sum_{k=1}^\infty \left(\dot{\eta}_i,\int_J g(\tau_k,x)(s)\left\langle c_k, \alpha_k(s)\varphi(\tau_k)(s)\right\rangle ds\right) \leq 0
\end{equation}
for any $\varphi \in \mathcal T$, $\varphi \geq 0$. Let $k_0 \in \mathbb N$ be fixed. Let $\varphi \in \mathcal T$ be such that $\hat{\varphi} \equiv 0$, $\varphi(\tau_k)(\cdot) \equiv 0$ ($k \ne k_0$), $\varphi(\tau_{k_0})(\cdot) \geq 0$. Then (\ref{vineq}) implies that
\begin{equation}
\label{vineq2}
\left(\dot{\eta}_i(x),\int_J g(\tau_{k_0},x)(s)\left\langle c_{k_0}, \alpha_{k_0}(s)\varphi(\tau_{k_0})(s)\right\rangle ds\right) \leq 0.
\end{equation}
Due to the linearity of the integral, the inequality (\ref{vineq2}) is equivalent to
\begin{equation}
\label{vs}
\int_J\bigl(\dot{\eta}_i(x),g(\tau_{k_0},x)(s)\langle c_{k_0},\alpha_{k_0}(s)\rangle\bigr)\varphi(\tau_{k_0})(s)ds \leq 0
\end{equation}
for any $\varphi(\tau_{k_0})(\cdot) \geq 0$. Then due to the continuity of the functions in (\ref{vs}) we have
\begin{equation*}
\bigl(\dot{\eta}_i(x),g(\tau_{k_0},x)(s)\langle c_{k_0},\alpha_{k_0}(s)\rangle\bigr) \leq 0
\end{equation*}
for all $s \in J$. Since $k_0 \in \mathbb N$ was chosen arbitrarily, we obtain (\ref{v2i}).

Let $\varphi \in \mathcal T$, $\varphi \geq 0$. Since the change of the dynamic values of $\varphi$ in finitely many points $\tau_k$ does not change $\hat{\varphi}=\ordin(\varphi)$, and the series in the right-hand side of (\ref{vineq}) converges, we obtain that
\begin{equation*}
\int_I\bigl(\dot{\eta}_i(x),f(t,x)+\hat{g}(t,x)w(t)\bigr)\hat{\varphi}(t)dt \leq 0
\end{equation*}
for all $\varphi \in \mathcal T$, $\varphi \geq 0$, which implies (\ref{v1i}). 

Since $x \in \partial M$, $i \in L_x$ were chosen arbitrarily, according to Theorem \ref{nagteo} we obtain that the conditions (\ref{v1i}) and (\ref{v2i}) imply that $M$ has the property of viability for
\begin{equation}
\label{veq5}
\dot{x}=f(t,x)+\hat{g}(t,x)w(t),
\end{equation}
and for the system (\ref{limitsystem}) (for any $k \in \mathbb N$), respectively.

2) Consider first the particular case where there exists $l>0$ such that $\tau_{k+1}-\tau_k \geq l$ for all $k \in \mathbb N$.
Let $\Omega=(0,\infty)$.

Let $\tau_1>0$. Then $x(0)=x_0$ and by Theorem \ref{teo1} since $\tau_{k+1}-\tau_k \geq l>0$ for all $k \in \mathbb N$ there exists
$\eta>0$ such that $x$ has the ordinary values $x(t)$ and
\begin{equation}
\label{eqabove}
\dot{x}(t)=f\bigl(t,x(t)\bigr)+\hat{g}\bigl(t,x(t)\bigr)w(t),
\end{equation}
for all $t \in [0,\eta)$. As is mentioned above, $M$ has the property of viability (\ref{eqabove}) on $(0,\eta)$, so
$x(t) \in M$ for all $t \in [0,\eta)$.

Now suppose that $\tau_1=0$. Then by Theorem \ref{teo1} the jump of the solution $x$ at $\tau_1=0$
can be found from (\ref{limitsystem}) at point $\tau_1$. Due to the remark above we have that $M$ has the property of viability for (\ref{limitsystem}) at the point $\tau_1$, so since $\gamma_1(-1/2)=x(0-) \in M$ we have that $\gamma_1(s) \in M$ for all $s \in J$.
Consequently, $x(0+)=\gamma(1/2) \in M$. Analogously to the first case we obtain that there exists 
$\eta>0$ such that $x(t)(s) \in M$ for all $t \in [0,\eta)$, $s \in J$.

We show that $x(t)(s) \in M$ for all $t \geq 0$, $s\in J$. By the change of the independent variable $t$ we obtain that the inclusion $x(t_0-) \in M$ implies that there exists $\eta=\eta(t_0,x(t_0-))>t_0$ such that
\begin{equation}
\label{viabincl}
x(t)(s) \in M 
\end{equation}
for all $t \in [t_0,\eta)$, $s \in J$. 
Thus, we obtain a strongly monotonically increasing sequence
$\{t_k\}_{k=1}^\infty$ such that $M$ has the property of viability on $(0,t_k)$
Suppose that the sequence $\{t_k\}_{k=1}^\infty$ is bounded from above. Consequently, $t_k \to t^*$ from the left, where $0<t^*<\infty$.
By our assumption the solution $x$ is defined for all $t \geq 0$. Since $M$ is closed, the limit $x(t^*-) \in M$. Thus, we may change the independent variable $t$, and apply the same arguments for $t_0=t^*$.
As a result, we obtain a contradiction with the assumption that $t^*$ is the maximal possible, so (\ref{viabincl}) holds for all $t \in (0,\infty)$, $s \in J$. 

The case of bounded $\Omega$ is treated similarly.

3) Consider the general case. Without loss of generality we give a proof for the case $\Omega=(0,\infty)$. Let $[c,d] \subset \Omega$, $u \in s\mathbb{BV}_n(c,d)$, $v=\dot{u} \in \mathcal T_n'(c,d)$, 
\begin{equation*}
u=q+\sum_{k=1}^\infty \langle c_k,\theta_{\tau_k}^{\beta_k}\rangle.
\end{equation*}
where $q:I \to \mathbb{R}^n$ is continuously differentiable, $w=\dot{q}$, $\beta_k \in \mathbb{AC}_n(J)$,
$\alpha_k=\dot{\beta}_k$ ($k \in \mathbb N$),
and without loss of generality $\tau_k \in (c,d)$ ($k \in \mathbb N$). Let us define
\begin{equation*}
u_l=q+\sum_{|c_k| \geq 1/l}\langle c_k,\theta_{\tau_k}^{\beta_k}\rangle,
\end{equation*}
where $v_l=\dot{u}_l \in \mathcal T_n'(c,d)$ contains the linear combination of delta-functions.
Since
\begin{equation*}
\|u-u_l\|_{s\mathbb{BV}_n(c,d)}=\biggl\|\sum_{|c_k| \leq 1/l} \langle c_k,\theta_{\tau_k}^{\beta_k}\rangle\biggr\|_{s\mathbb{BV}_n(c,d)} \to 0
\end{equation*}
($l \to \infty$), by Theorem \ref{existsteo2} we have the convergence $x_l \to x$
in $s\mathbb{BV}_n(c,d)$ of the sequence of solutions $x_l \in s\mathbb{BV}_n(c,d)$ of the initial value problems (\ref{eq11})
for $v_l=\dot{u}_l$. We apply the results obtained above to the initial value problems (\ref{eq11})
corresponding to $v_l \in \mathcal T_n'$, so, any solution $x_l$ is viable in $M$ on $(c,d)$, i.e., $x_l(t)(s) \in M$ for all $t \in (c,d)$, $s \in J$. Convergence $x_l \to x$
in $s\mathbb{BV}_n(c,d)$ implies that
\begin{equation*}
x_l(t)(s) \to x(t)(s) \quad \bigl(t \in (c,d),~s \in J\bigr).
\end{equation*}
Since $M$ is closed, we have $x(t)(s) \in M$ for all $t \in (c,d)$, $s \in J$. Since
$[c,d] \subset \Omega$ was chosen arbitrarily, we obtain that $x(t)(s) \in M$ for all $t \geq t_0$, $s \in J$.
\end{proof}

\begin{corollary}
\label{viab_cor}
Let $\dot{\eta}_i(x) \in \mathbb{R}^n$ \rm(\it $i \in L_x$\rm) \it be linearly independent for any $x \in \partial M$.
If
\begin{equation}
\label{v1ia}
\bigl(\dot{\eta}_i(x),f(t,x)+\hat{g}(t,x)w(t)\bigr)\leq 0
\end{equation}
for all $t \in \Omega$, $x \in \partial M$,
\begin{equation}
\label{v2ib}
\bigl(\dot{\eta}_i(x),g(\tau_k,x)(s)\bigl\langle c_k, \alpha_k(s)\bigr\rangle\bigr)\leq 0
\end{equation}
for all $s \in J$, $k \in \mathbb N$, $x \in \partial M$, then
$M$ has the property of viability for \rm (\ref{eq11}) \rm(\it on $\Omega$\rm).
\end{corollary}
\begin{proof}
The proof follows from the proof of Theorem \ref{viabteo}.
\end{proof}

\begin{example} Let $I=(-1,\infty)$, $D=I \times (-1,2)$. Let $\Omega=(0,\infty) \subset I$, $\{\tau_k\}_{k=1}^\infty \subset \Omega$,
$\tau_k \to \infty$. Let us consider in $\mathcal T'$ the following ordinary differential equation with distributions of the form (\ref{eq11}),
\begin{equation}
\label{eqexviab}
\dot{x}=-x+\biggl(\frac{1}{2}-x\biggr)\sum_{k=1}^\infty \delta_{\tau_k}^{\alpha_k}, \quad \alpha_k \geq 0, \quad \alpha_k \in \mathbb C(J).
\end{equation}
Let $\eta(x)=(x-1/2)^2-1/4$, so 
\begin{equation*}
M=\{x \in \mathbb R: \eta(x) \leq 0\}=[0,1]. 
\end{equation*}
Notice, that $\dot{\eta}(x) \in \{-1, 1\}$ if $x \in \partial M=\{0,1\}$. Since the system (\ref{eqexviab}) has form (\ref{eq2}) for $f(x)=-x$, $g(x)=1/2-x$, and 
$v=\sum_{k=1}^\infty \delta_{\tau_k}^{\alpha_k}$,
we may apply Corollary \ref{viab_cor}. We have
\begin{equation*}
\dot{\eta}(1)f(1) \leq 0, \quad \dot{\eta}(0)f(0) \leq 0,
\end{equation*}
and
\begin{equation*}
\dot{\eta}(1)g(1)\alpha_k(s) \leq 0, \quad \dot{\eta}(0)g(0)\alpha_k(s) \leq 0 \\ [2mm]
\end{equation*}
for all $s \in J$, $k \in \mathbb N$, where $\dot{\eta}(1)=1$, $g(1)=-1/2$, $\dot{\eta}(0)=-1$, $g(1)=-1/2$, $f(1)=-1$, $f(0)=0$ and $\alpha_k \geq 0$ ($s \in J$, $k \in \mathbb N$). According to Corollary \ref{viab_cor} the set $M$ has the property of viability for the system (\ref{eqexviab}) on $\Omega$. 
\end{example}

%

Let us consider the applications of Theorem \ref{viabteo} to stability analysis. 

Let $I=(a,\infty)$. Consider the ordinary differential equation with distributions of the form (\ref{eq1}),
\begin{equation}
\label{stabeq1}
\dot{x}=f(x)+g(x)v,
\end{equation}
where $f:D \to \mathbb R^n$, $g:D \to \mathbb R^{n \times n}$ are Lipschitz in $D$. The solution of (\ref{stabeq1}) which is identically equal to a constant $x^* \in D$ is called the \textit{equilibrium point} (clearly, $x^*$ is an equilibrium point of (\ref{stabeq1}) if and only if $f(x^*)=g(x^*)=0$). 

Following the standard terminology, we note that the equilibrium point $x^* \in D$ is \textit{uniformly stable}, if there exists a sequence of the closed subsets $M_{n+1} \subset M_k \subset D$ ($k \in \mathbb N$), where
$M_k$ contains an open neighbourhood of $x^* \in D$ and is also contained in another open neighbourhood of $x^* \in D$, such that for any $l \in \mathbb N$ there exists $k \geqslant l$ having the property: any solution $x$ of (\ref{stabeq1}) with $x(t_0-) \in M_k$ is viable in $M_l$ on $\Omega=(t_0,\infty) \subset I$.

\begin{theorem}
\label{stabteo}
Suppose that $x^* \in D$ is an equilibrium point of {\rm(\ref{stabeq1})}, and
\begin{equation}
\label{stabineq}
\bigl(x-x^*,f(x)+g(x)v\bigr) \leqslant 0
\end{equation}
in $\mathcal T_n'(\Omega)$ for any $x \in D$ such that $|x-x^*|_2=1/l$, where $|\cdot|_2$ is a Euclidean norm in $\mathbb R^n$, $l \in \mathbb N$. Then $x^*$ is uniformly stable.
\end{theorem}
\begin{proof}
Let us note that in the definition of the uniform stability above we may have $l=k$, that is, it suffices to prove that for each $l \in \mathbb N$ the subset $M_l$ has the property of viability for (\ref{viabeq1}) on $\Omega=(t_0,\infty)$. Let
\begin{equation*}
M_l=\{x \in D:\beta_l(x) \leqslant 0\}, \text{ where } \beta_l(x)=|x-x^*|_2^2-\frac{1}{l^2}.
\end{equation*}
Then the sequence $\{M_l\}_{l=1}^\infty$ satisfies the conditions above. We have
\begin{equation*}
\dot{\beta}_l(x)=2(x-x^*) \ne 0
\end{equation*}
for any $x \in D$ such that $|x-x^*|_2=1/l$. According to Theorem \ref{viabteo} the inequality (\ref{stabineq}) implies that $M_l$ possesses the property of viability for (\ref{stabeq1}) on $\Omega$, so $x^*$ is uniformly stable.
\end{proof}

Analogously to Theorem \ref{viabteo}, we obtain the following corollary of Theorem \ref{stabteo}.

\begin{corollary}
\label{stabcor}
Suppose that $x^* \in D$ is an equilibrium point of {\rm(\ref{stabeq1})}, and
\begin{equation}
\label{stabineq2}
\bigl(x-x^*,f(x)+g(x)w(t) \bigr) \leqslant 0
\end{equation}
for all $t \in \Omega=(t_0\infty)$ and all $x \in D$ such that $|x-x^*|_2=1/l$, 
\begin{equation}
\label{stabineq3}
\bigl(x-x^*,g(\tau_k,x)(s)\bigl\langle c_k, \alpha_k(s)\bigr\rangle\bigr) \leqslant 0
\end{equation}
for all $s \in J$, $k \in \mathbb N$ and all $x \in D$ such that $|x-x^*|_2=1/l$, 
where $|\cdot|_2$ is a Euclidean norm in $\mathbb R^n$, $l \in \mathbb N$. Then $x^*$ is uniformly stable.
\end{corollary}
\begin{proof}
The proof follows from the proof of Theorem \ref{viabteo}.
\end{proof}

\begin{example}
Let us consider the ordinary differential equation with distributions
\begin{equation}
\label{stabeq2}
\dot{x}=\biggl(\frac{1}{2}-x\biggr)\sum_{k=1}^\infty \delta_{\tau_k}^{\alpha_k},
\end{equation}
where $I=(-1,\infty)$, the shape $\alpha_k \geqslant 0$ is continuous. Clearly, (\ref{stabeq2}) has an equilibrium point
\begin{equation*}
x_*=\frac{1}{2}. 
\end{equation*}
The equation (\ref{stabeq2}) has the form of (\ref{stabeq1}) for
$f(x) \equiv 0$, $g(x)=\frac{1}{2}-x$, $w(t) \equiv 0$. 
We have to show that the inequality (\ref{stabineq3}) holds. Indeed, the inequality (\ref{stabineq3}) is equivalent to
\begin{equation*}
\left(\frac{1}{2}-x\right)\alpha_k(s)\left(x-\frac{1}{2}\right)=-\left(\frac{1}{2}-x\right)^2 \alpha_k(s) \leqslant 0 
\end{equation*}
for any $s \in J$ and any $x$ such that $|x-1/2|=1/l$ ($l \in \mathbb N$). According to Theorem \ref{stabteo} the equilibrium point $x_*=\frac{1}{2}$ of equation (\ref{stabeq2}) is uniformly stable.
\end{example}

\section{The statement of the impulse problem of avoidance of encounters}

Let $M \subset \mathbb R^n$ be a closed subset, $D=I \times N$, where $N$ is an open subset,
$M \subset N$.
Let us consider in $\mathbb R^n$ the following controlled system of the form (\ref{intro_contr2}),
\begin{equation}
\label{control2}
\dot{x}=f(t,x)+g(t,x)v, \quad x(t_0)=x_0, \quad v \in \mathcal V, \\ [2mm]
\end{equation}
where $f$, $g$ satisfy conditions 1), 2), the function $f$ is continuous in $t$, and the set of admissible ordinary controls $\mathcal V$ is given by
\begin{equation}
\label{v1}
\mathcal V=\left\{v \in \mathbb L_n(I): v_i(\cdot) \geq 0, \int_{\Omega} v_i(s)ds \leq V, 1 \leq i \leq n\right\},
\end{equation}
where $\Omega=(t_0,T) \subset I$, $V>0$ is given, $v=(v_i)_{i=1}^n$.

Let $\mathcal V_M(T) \subset \mathcal V$ be the maximal set of admissible controls such that
\begin{equation*}
v \in \mathcal V_M(T) \text{ implies that } x \text{ is viable in } M \text{ on } \Omega=(t_0,T) \subset I,
\end{equation*}
where $x \in \mathbb{AC}_n$ is the solution of system (\ref{control2}).
According to \cite{Faz}, we call the following maximization problem
\begin{equation}
\label{max2}
T \underset{v}{\to} \max, \quad
v \in \mathcal V_M(T)
\end{equation}
the \textit{problem of avoidance of encounters with the set }$\mathbb{R}^n \setminus M$. 

Along with the system (\ref{control2}), let us consider in $\mathcal T_n'$ the following controlled system 
\begin{equation}
\label{control3}
\dot{x}=f(t,x)+g(t,x)v, \quad
x(t_0-)=x_0, \quad
v \in \mathcal V',
\end{equation}
where $f$, $g$ satisfy 1), 2), the function $f$ is continuous in $t$, and the set of admissible distributional (i.e., \textit{impulse}) controls $\mathcal V'$ is given by
\begin{equation}
\label{v2}
\mathcal V'=\left\{v \in \mathcal T_n'(I) : v_i \geq 0, \int_{t_0}^{T} v_ids \leq V, 1 \leq i \leq n\right\},
\end{equation}
where $V>0$ is given, $v=(v_i)_{i=1}^n$ (the definitions of a non-negative distribution and integral of a distribution were given in Section 2).
Notice that $\mathcal V \subset \mathcal V'$, where the elements of $\mathcal V$ are viewed as the regular distributions in $\mathcal T_n'$. Thus, system (\ref{control3})(\ref{v2}) is an extension of system (\ref{control2})(\ref{v1}). 

We define $\mathcal V_M'(T) \subset \mathcal V'$ to be the maximal set of admissible controls such that 
\begin{equation*}
v \in \mathcal V'_M(T) \text{ implies that } x \text{ is viable in } M \text{ on } \Omega=(t_0,T) \subset I,
\end{equation*}
where $x \in s\mathbb{BV}_n$ is a solution of system (\ref{control3}).
Analogously to \cite{Faz}, we call
the maximization problem
\begin{equation}
\label{max2}
T \underset{v}{\to} \max, \quad
v \in \mathcal V'_M(T)
\end{equation}
the \textit{impulse problem of avoidance of encounters with the set }$\mathbb{R}^n \setminus M$. \\

As follows from the next example, the problem of avoidance of encounters may have no solution $(T^*,v^*)$ for the system (\ref{control2})(\ref{v1}), but may have a solution for the extended system (\ref{control3})(\ref{v2}).

\begin{example}
Let $I=(-1,\infty)$, $D=I \times \mathbb R$. Let $\eta(x)=x^2-1$, so 
\begin{equation*}
M=\{x \in \mathbb R: \eta(x) \leq 0\}=[-1,1]. 
\end{equation*}
Let us consider in $\mathcal T'$ the following controlled system
\begin{equation}
\label{control4}
\dot{x}=x-v, \quad v \in \mathcal V, \quad x(0-)=1,
\end{equation}
\begin{equation*}
\mathcal V=\left\{v \in \mathcal T':v \geq 0, \int_0^1v dt \leq \frac{1}{2}\right\}.
\end{equation*}
Let us show that the solution of the problem of avoidance of encounters with the set $\mathbb R \setminus M$ for (\ref{control4}) is given by $\Omega^*=(0,T^*) \subset I$,
\begin{equation*}
T^*=\ln(2), \quad v^*=\frac{1}{2}\delta_0^\alpha, \quad \alpha \geq 0.
\end{equation*}
Suppose that $v=1/2\delta_0^\alpha$, $\alpha \geq 0$. Then according to Theorem \ref{teo1} we have $x(0)(s)=\gamma(s)$ ($s \in J$), $x(t)=1/2e^{t}$, $t \in \bigl(0,\ln(2)\bigr)$, where
\begin{equation*}
0 \leq \gamma(s)=1-\frac{1}{2}\int_{-\frac{1}{2}}^s\alpha(\eta)d\eta \leq 1
\end{equation*}
for all $s \in J$, so $x$ is viable in $M$ on $\Omega=(0,\ln(2)) \subset I$. 
Along with that, for any regular control $v$ solution of (\ref{control4}) is given by
\begin{equation*}
x(t)=e^t-\frac{e^t}{2}\int_0^t v(s)e^{-s}ds
\end{equation*}
for all $t>0$.
As follows from the obtained representation, there exists $\xi>0$ such that 
$x(t) \geq e^t-e^t\frac{e^{-\xi}}{2}$ 
for all $t \geq 0$. Thus, 
\begin{equation*}
T<\ln\biggl(\frac{1}{1-\frac{e^{-\xi}}{2}}\biggr)<\ln(2), 
\end{equation*}
so $T^*=\ln(2)$ is the maximal viability time.
\end{example}

In the subsequent paper we provide the necessary conditions for optimality in the \textit{impulse problem of avoidance of encounters with the set} $\mathbb R^n \setminus M$.

\bibliographystyle{plain}
\bibliography{eng_article_diffeq}

\end{document}